\documentclass[twoside,leqno]{article}

\usepackage{amssymb}
\usepackage{amsmath}
\usepackage[cp1250]{inputenc}
\usepackage{enumerate}

\numberwithin{equation}{section}
\newtheorem{thm}[equation]{\sc Theorem}
\newtheorem{lem}[equation]{\sc Lemma}

\newtheorem{rem}[equation]{\sc Remark }

\newtheorem{ex}[equation]{\sc Example}

\makeatletter
\renewcommand{\@seccntformat }[1]{\csname the#1\endcsname. }
\makeatother

 1

\def\CH{{\cal H}}

\def\CM{{\cal M}}

\def\CP{{\cal P}}
\def\CO{{\cal O}}

\def\vsp{\vspace*{1.5ex}}

\def\epv {{$\mbox{}$\hfill ${\Box}$\vspace*{1.5ex} }}

\def\longarr#1#2{{\buildrel{#1} \over {\hbox to #2pt{\rightarrowfill}}}}

\def\id{\mbox{{\rm id}}}

\def\mapup#1{\Big\uparrow\rlap{$\vcenter{\hbox{$\scriptstyle#1$}}$}}

\def\End{\mbox{\rm End}}

\def\Gl{\mbox{\rm Gl}}

\def\vsp{\vspace*{1.5ex}}

\def\vedge#1{{\buildrel{#1} \over {\hbox to
20pt{\hspace{-0.2em}$-$\hspace{-0.2em}$-$\hspace{-0.2em}$-$ }}}}

\begin{document}

\renewcommand{\thefootnote}{}
\def\refname{\begin{center}\normalsize{\it REFERENCES}\end{center}}

 \title{\large{\bf Generic extensions of nilpotent $k[T]$-modules, monoids of partitions and constant terms of Hall polynomials}
   \\
   \vspace{0.5em}}
\author{  \normalsize {\bf Justyna  Kosakowska}
\\
 \small{\it Faculty of Mathematics and Computer Science,
Nicolaus Copernicus University}\\ \small{\it ul. Chopina 12/18,
87-100 Toru\'n, Poland} \\
\small{\it E-mail: justus@mat.uni.torun.pl}}
\date{}
\maketitle


\footnote{Partially supported by Research Grant No. N N201 269135 of Polish Ministry of Science and High Education
  and by Grant 403-M of UMK}

\footnotesize {\bf Abstract.} We prove that the monoid of generic extensions of finite dimensional nilpotent $k[T]$-modules is isomorphic to
the monoid of partitions (with addition of partitions). Moreover we give a~combinatorial algorithm that calculates constant terms
of classical Hall polynomials.   \vsp \vsp



\normalsize

\section{Introduction}

Let $k$ be an~algebraically closed field and let $k[T]$ be the $k$-algebra of polynomials
in one variable $T$. We consider nilpotent $k[T]$-modules $M$, $N$  and the generic extension $M\ast N$ of $M$ by $N$, i.e.
an~extension of $M$ by $N$ with the minimal dimension of its
endomorphism ring (see Section \ref{sec-def} for definitions). By results presented in
\cite{bon,dengdu,reineke} generic extensions of nilpotent $k[T]$-modules exist and the operation of taking the generic
extension provides the set of all isomorphism classes of
nilpotent $k[T]$-modules with a~monoid structure $\CM^\ast$. In this paper we study connections of this
monoid with the monoid $\CP^+$ of all partitions with addition of partitions as an action. More precisely,
we prove in Theorem \ref{thm-main} that these monoids are isomorphic. This isomorphism gives us
a~combinatorial description of generic extensions that have a~geometric nature. For 
a~geometric interpretation of generic extensions the reader is referred to \cite{reineke,reineke1}.

On the other hand there is
a~$\mathbb{C}$-algebra isomorphism $\mathbb{C}\CM^\ast\simeq \CH_0$, where 
$\CH_0$ is the specialisation of the
Hall algebra $\CH_q$ to $q=0$ and $\mathbb{C}\CM^\ast$ is the $\mathbb{C}$-algebra generated
by the monoid $\CM^\ast$ (see \cite{dengdu,hubery,wolf} and Section \ref{sec-monoids}).
There are many results that show connections between generic extensions, Hall polynomials and Ringel-Hall algebras
(see \cite{reineke, reineke1, hubery, wolf}  for Dynkin and extended Dynkin quivers, \cite{dengdu} for
cyclic quivers, \cite{kos04,kos05a} for poset representations).

In Section \ref{sec-algorithm}, exploring the isomorphism $\mathbb{C}\CM^\ast\simeq \CH_0$ (explicitly given in \cite{wolf}),
we describe a~combinatorial algorithm that finds the constant terms of classical Hall polynomials.

{\bf Acknowledgments.} The author would like to thank Stanis\l aw Kasjan for his helpful  
comments concerning the proof of Lemma \ref{lem-ext}. 

\section{Notation and definitions} \label{sec-def}

Throughout this paper $k$ is a~fixed algebraically closed field.

Let $\lambda=(\lambda_1,\ldots,\lambda_n,\ldots)$ be a~partition (i.e. a~sequence  of
non-negative integers containing only finitely many
non-zero terms and such that $\lambda_1\geq\lambda_2\geq\ldots\geq\lambda_n\geq\ldots$). Denote by
$\overline{\lambda}=(\overline{\lambda}_1,\ldots,\overline{\lambda}_n,\ldots)$ the dual partition
of $\lambda$,
i.e.  $$\overline{\lambda}_i=\#\{j\; ;\; \lambda_j\geq i\},$$ where $\#X$ denotes the cardinality of a~finite set $X$.
We identify partitions that differ only by a~string of zeros at the end. 
Let $\CP$ be the set of all partitions. Denote by $|\lambda|$ the \textbf{weight} of $\lambda$ defined by 
$$|\lambda|=\lambda_1+\lambda_2+\ldots$$
and by $0=(0)$ the unique partition of zero. Consider two associative monoids:
\begin{itemize}
 \item $\CP^+=(\CP,+,0)$, where 
$(\lambda_1,\lambda_2,\ldots)+(\nu_1,\nu_2,\ldots)=(\lambda_1+\nu_1, \lambda_2+\nu_2,\ldots)$;
\item $\CP^\cup=(\CP,\cup,0)$, where 
$(\lambda_1,\lambda_2,\ldots)\cup(\nu_1,\nu_2,\ldots)=(\mu_1,\mu_2,\ldots)$ and
$(\mu_1,\mu_2,\ldots)$ is the partition that is consisted of integers $\lambda_1,\lambda_2,\ldots,\nu_1,\nu_2,\ldots$
arranged in the descending order (e.g. $(3,3,2,1)+(2,2)=(5,5,2,1)$ and $(3,3,2,1)\cup(2,2)=(3,3,2,2,2,1)$);
\end{itemize}
By \cite[1.8]{macd} the operations $+$ and $\cup$ are dual to each other 
(i.e. $\overline{\lambda\cup\nu}=\overline{\lambda}+\overline{\nu}$). 
One of the main aims of the paper is to describe connections of these monoids with 
extensions of nilpotent $k[T]$-modules.

Let $k[T]$ be the $k$-algebra
of polynomials in the variable $T$. For any partition $\lambda=(\lambda_1,\ldots,\lambda_n,\ldots)$,
where $\lambda_{n+1}=\lambda_{n+2}=\ldots=0$,
denote by $$M(\lambda)=M(\lambda,k)\cong k[T]/(T^{\lambda_1})\oplus\ldots\oplus k[T]/(T^{\lambda_n})$$
the corresponding $k[T]$-module. It is obvious that the function $\lambda\to M(\lambda)$
gives a~bijection between the set $\CP$ of all partitions and the set of all isomorphism classes of 
nilpotent $k[T]$-modules (i.e. finitely generated $k[T]$-modules $M$ such that $T^aM=0$
for some $a\geq 0$). Denote by $\CM$ a~set of representatives of all isomorphism classes of nilpotent $k[T]$-modules.

Let $M,N\in \CM$. By \cite{bon}, \cite{dengdu} and \cite{reineke}, there is the~unique (up to isomorphism) extension $X$ of $M$ by $N$
with the minimal dimension of endomorphism ring $\End_{k[T]}(X)$, i.e. 
a~nilpotent $k[T]$-module $X$ such that there exists a~short exact sequence of the form
$$ 0\to N\to X\to M\to 0.$$ The module $X$ is called the \textbf{generic extension} of $M$ by $N$ and
is denoted by $X=M\ast N$. Denote by $M\oplus N$ the direct sum of the modules $M$ and $N$ and by $0$
the unique zero module.
Consider two monoids:
\begin{itemize}
 \item $\CM^\ast=(\CM,\ast,0)$ (\textbf{the monoid of generic extensions}),
 \item $\CM^\oplus=(\CM,\oplus,0)$.
\end{itemize}
The associativity of the monoid $\CM^\ast=(\CM,\ast,0)$ follows by \cite{dengdu},
whereas that of the monoid $\CM^\oplus=(\CM,\oplus,0)$
is obvious.

\section{Generic extensions and partitions}\label{sec-monoids}

The following fact is one of the main results of the paper.

\begin{thm}
 The function $$\Phi:\CP\to \CM $$
  such that $\Phi(\lambda)=M(\lambda)$, for any partition $\lambda$, 
induces isomorphisms of monoids:
  $$ \Phi:\CP^+\to \CM^\ast $$
and 
  $$ \Phi:\CP^\cup\to \CM^\oplus .$$ 
Moreover 
$$ M(\overline{\alpha})\ast M(\overline{\beta})=M(\overline{\alpha\cup\beta}).$$
\label{thm-main}\end{thm}

To prove Theorem \ref{thm-main} we need a~geometric interpretation of generic extensions.

We  identify $k[T]$-modules of the form  $M(\lambda,k)$ with systems $M(\lambda,k)=(V,\varphi)$,
where $V$ is a~finite dimensional $k$-vector space and $\varphi:V\to V$ is a~nilpotent linear
endomorphism with the Jordan type $\lambda$ (i.e. nilpotent representation of a~loop quiver). 
By $\mathcal{N}(k)$ we denote the category of all such systems.
If $(V,\varphi)$, $(V_1,\varphi_1)$ are objects in $\mathcal{N}(k)$, then a~morphism 
$f:(V,\varphi)\to (V_1,\varphi_1)$ is a~linear map $f:V\to V$ such that $\varphi_1f=f\varphi$.
It is easy to see that the category $\mathcal{N}(k)$ is equivalent to the category of all finite dimensional nilpotent
$k[T]$-modules. For an~account of the theory of modules and quiver representations we refer the reader to \cite{ASS} and \cite{ARS}. 

Consider the affine variety $\mathbb{M}_n(k)$ of all $n\times n$-matrices with coefficients
in $k$. The general linear group $\Gl_n(k)$ acts on $\mathbb{M}_n(k)$ via conjugations, i.e. for $g\in\Gl_n(k)$
and $M\in \mathbb{M}_n(k)$, we put $g\cdot M=gMg^{-1}$. Let $\mathbb{M}_n^{nil}(k)$ be the~subset of $\mathbb{M}_n(k)$
consisted of all nilpotent matrices. The subset $\mathbb{M}_n^{nil}(k)$ is closed in $\mathbb{M}_n(k)$ 
(in Zariski topology) and it is closed under the action of  $\Gl_n(k)$. It is easy to observe that
points of $\mathbb{M}_n^{nil}(k)$ corresponds bijectively to the objects  $(V,\varphi)$ of $\mathcal{N}(k)$
with $\dim_kV=n$. Moreover the orbits of the action of $\Gl_n(k)$ on $\mathbb{M}_n^{nil}(k)$
corresponds bijectively to the isomorphism classes of the objects $V$ in $\mathcal{N}(k)$ (with $\dim_kV=n$) and hence to 
the isomorphism classes of finite dimensional nilpotent $k[T]$-modules 
$V$ (with $\dim_kV=n$). If $M(\lambda)\equiv(V,\varphi)$ is a~nilpotent $k[T]$-module
with $\dim_kM(\lambda)=n$,
then we denote by $\CO_\lambda$ (resp. $\overline{\CO_\lambda}$) the orbit (resp. the Zariski-closure)
 of $\varphi\in \mathbb{M}_n^{nil}(k)$ of the $\Gl_n(k)$-action.  
 
Let $\lambda,\nu$ be partitions with weights $|\lambda|=|\nu|=n$. We say that a~module $M(\lambda)$ \textbf{degenerates} to
the module $M(\nu)$, if $\CO_{\nu}\in \overline{\CO_{\lambda}}$. If $M(\lambda)$ degenerates to
$M(\nu)$ we write $M(\lambda)\leq_{deg} M(\nu)$. The relation $\leq_{deg}$ is a~partial order
on isomorphism classes of finite dimensional nilpotent $k[T]$-modules.  Geometrically, the generic extension
$M(\lambda)\ast M(\nu)$ (resp. the direct sum $M(\lambda)\oplus M(\nu)$) is the 
$\leq_{deg}$-minimal (resp. $\leq_{deg}$-maximal) extension of $M(\nu)$ by $M(\mu)$, i.e. if $X$ is
an~extension of $M(\nu)$ by $M(\mu)$, then $M(\lambda)\ast M(\nu)\leq_{deg} X$ 
(resp. $X\leq_{deg}M(\lambda)\oplus M(\nu)$), see \cite{bon}, \cite{dengdu} and \cite{reineke}.
For an~introduction to geometric methods in representation
theory the reader is referred to \cite{kraft1} and \cite{bon}.

The following fact is proved in \cite[I.3]{kraft}.

\begin{thm}
Let $\lambda,\nu$ be partitions with $|\lambda|=|\nu|$. $M(\lambda)\leq_{deg}M(\nu)$ if and only if, for any $m\geq 1$:
$$ \sum_{i=1}^m\overline{\lambda}_i \leq \sum_{i=1}^m\overline{\nu}_i.$$
\label{thm-deg}\end{thm}

The following lemma is used in the proof of Theorem \ref{thm-main}.

\begin{lem}
 Let $\sigma,\nu,\mu$ be partitions. If there exists a~short exact sequence
$$ 0\longarr{}{30}M(\nu)\longarr{a}{30}M(\sigma)\longarr{b}{30}M(\mu)\longarr{}{30}0,$$
then for any $m\geq 1$:
  $$\sum_{i=1}^m\sigma_i\leq\sum_{i=1}^m\lambda_i, $$
where $\lambda=\mu+\nu$.
\label{lem-ext}\end{lem}

 \textbf{Proof.} The proof is by induction on $|\nu|$. If $|\nu|=0$,
then $M(\sigma)\cong M(\mu)$, $\sigma=\mu$ and we are done.

Assume that $|\nu|>0$. We have $\nu=(\nu_1,\ldots,\nu_n)$, $\nu_n\neq 0$ and
$$  M(\nu)=M(\nu_1)\oplus\ldots\oplus M(\nu_n).$$
Consider the monomorphism
  $$ f=[\iota,0\ldots,0]:M(1)\to M(\nu_1)\oplus\ldots\oplus M(\nu_n),$$
where $\iota:M(1)\to M(\nu_1)$ is an~inclusion. By the Snake Lemma, we get the following diagram with
exact rows and columns:
$$ 
\begin{array}{ccccccccc}
  &&0&&0&&0&& \\
  &&\mapup{} && \mapup{}&&\mapup{} && \\
  0&\longarr{}{30}&M(\nu')&\longarr{}{30}&M(\sigma')&\longarr{}{30}&M(\mu)&\longarr{}{30}&0\\
  &&\mapup{} && \mapup{}&&\mapup{} && \\
  0&\longarr{}{30}&M(\nu)&\longarr{a}{30}&M(\sigma)&\longarr{b}{30}&M(\mu)&\longarr{}{30}&0\\
  &&\mapup{f} && \mapup{a\cdot f}&&\mapup{0} && \\
  0&\longarr{}{30}&M(1)&\longarr{\id}{30}&M(1)&\longarr{0}{30}&0&\longarr{}{30}&0\\
  &&\mapup{} && \mapup{}&& && \\
  &&0&&0&&&&
\end{array}
$$ 
where $\nu'=(\nu_1-1,\nu_2,\ldots,\nu_{n-1},\nu_n)$ and there exists $i$ such that $\sigma'_i=\sigma_i-1$
and $\sigma'_j=\sigma_j$ for $j\neq i$. By the induction hypothesis we get
$$\sum_{i=1}^m\sigma_i'\leq\sum_{i=1}^m\lambda_i', $$
where $\lambda'=\mu+\nu'$. Therefore 
$$\sum_{i=1}^m\sigma_i\leq\sum_{i=1}^m\lambda_i, $$
where $\lambda=\mu+\nu$ and we are done.
\epv

\begin{lem}
 Let $\nu,\mu$ be partitions. We have
 $$ M(\nu)\ast M(\mu)=M(\nu+\mu), $$
where $ M(\nu)\ast M(\mu)$ is the generic extension of $M(\nu)$ by $M(\mu)$.
\label{lem-gen-ext}\end{lem}

\textbf{Proof.} It is easy to see that $M(\nu+\mu)$
is an~extension of $M(\nu)$ by $M(\mu)$. If $M(\sigma)$ is any extension of $M(\mu)$ by $M(\nu)$,
then by Lemma \ref{lem-ext}, for any $m\geq 1$:
  $$\sum_{i=1}^m\sigma_i\leq\sum_{i=1}^m\lambda_i, $$
where $\lambda=\mu+\nu$. By \cite[1.11]{macd}, for any 
$m\geq 1$:
  $$\sum_{i=1}^m\overline{\sigma}_i\geq\sum_{i=1}^m\overline{\lambda}_i. $$
Theorem \ref{thm-deg} yields $$M(\nu+\mu)=M(\lambda)\leq_{deg} M(\sigma). $$
Since $M(\lambda)\ast M(\nu)$ is the $\leq_{deg}$-minimal extension of $M(\nu)$ by $M(\mu)$, 
we are done.\epv

\textbf{Proof of Theorem \ref{thm-main}.} Let $$\Phi:\CP\to \CM $$
  be such that $\Phi(\lambda)=M(\lambda)$, for any partition $\lambda$.
By Lemma \ref{lem-gen-ext}, the induced function  
  $$ \Phi:\CP^+\to \CM^\ast $$
is an~isomorphism of monoids. It is easy to see that
  $$ \Phi:\CP^\cup\to \CM^\oplus .$$ 
is an~isomorphism of monoids.
Moreover $$ M(\overline{\alpha})\ast M(\overline{\beta})=M(\overline{\alpha}+\overline{\beta})=M(\overline{\alpha\cup\beta}),$$
because $\overline{\alpha}+\overline{\beta}=\overline{\alpha\cup\beta}$.
 \epv

\section{Constant terms of Hall polynomials}\label{sec-algorithm}

In this section we describe a~combinatorial algorithm that finds the constant
term of a~given Hall polynomial. 

Let $\alpha,\beta,\gamma$ be partitions and let $k$ be a~finite field. Denote by
  $$ F_{\alpha,\beta}^\gamma(k)$$
the number of submodules $U$ of $M(\gamma,k)$ that are isomorphic to $M(\beta,k)$
and the factor module $M(\gamma,k)/U$ is isomorphic to $M(\alpha,k)$. By the result of Hall (see \cite[II.4.3]{macd}),
there exists a~polynomial $\varphi_{\alpha\beta}^\gamma$ with integral coefficients such that:
  $$\varphi_{\alpha\beta}^\gamma(\# k) =F_{\alpha,\beta}^\gamma(k)$$
for any finite field $k$.
We call $\varphi_{\alpha,\beta}^\gamma$ the~\textbf{Hall polynomial} associated with partitions $\alpha,\beta,\gamma$.

By \cite{dengdu}, \cite{hubery} and \cite{wolf}, the complex algebra
$\mathbb{C}\CM^\ast$ generated by the monoid $\CM^\ast$ of generic extension
is isomorphic to the degenerate complex Hall algebra $\CH_0$, where  $\CH_0$ has a~basis 
$$\{u_{\alpha}\; ;\; \alpha\in\CP\} $$ as a~$\mathbb{C}$-vector
space
and the multiplication is given by the formula
  $$ u_\alpha u_\beta=\sum_\gamma \varphi_{\alpha\beta}^\gamma(0)u_\gamma. $$
  By \cite{wolf}, the isomorphism
  $$ F: \mathbb{C}\CM^\ast\to \CH_0$$
is given by the formula $$ F(M(\alpha))=\sum_{\beta\, :\;M(\alpha)\leq_{deg} M(\beta)}u_{\beta}. $$

We use the following notation. A~partition $\alpha=(\alpha_1,\ldots,\alpha_n,\ldots)$ shall be written as
  $$ (\ldots ,r^{m_r},\ldots ,2^{m_2}, 1^{m_1}), $$
where $m_r$ indicates the number of times the integer $r$ occurs in $\alpha$, e.g. $$(3,3,2,2,2,1,1,1,1)=(3^2,2^3,1^4).$$

\begin{lem}
 Let $\gamma$ be an~arbitrary partition and let $\alpha=(1^n),\beta=(1^m)$ be partitions  
with the property  $\varphi_{\alpha,\beta}^\gamma\neq 0$. Then
 $$ \varphi_{\alpha,\beta}^\gamma(0)=1. $$
\label{lem-constant-term}\end{lem}

\textbf{Proof.} 
Note that $F(M(\alpha))=u_\alpha$ and $F(M(\beta))=u_\beta$. Then
  $$ F(M(\alpha)) F(M(\beta))=u_\alpha u_\beta=\sum_\delta \varphi_{\alpha\beta}^\delta(0)u_\delta$$
  and $$  F(M(\alpha)) F(M(\beta))=F(M(\alpha)\ast M(\beta))=F(M(\alpha+\beta))= 
\sum_{M(\alpha+\beta)\leq_{deg} M(\delta)}u_{\delta}$$
Comparing these sums we get $\varphi_{\alpha,\beta}^\gamma(0)=1$, if $\varphi_{\alpha,\beta}^\gamma\neq 0$.
\epv

Applying recursively (following $\leq_{deg}$-order) methods applied in the proof of Lemma \ref{lem-constant-term} one can 
calculate constant terms
of Hall polynomials. We illustrate this algorithm in the following example.

\begin{ex}{\rm 
 We calculate the constant term of the Hall polynomial $\varphi_{(2,1)(2)}^{(4,1)}$. We  apply Theorem \ref{thm-deg}
and the definition of $F$.

\textbf{Step 1.} By Lemma \ref{lem-constant-term}, we have 
   $$\varphi_{(1^3)(1^2)}^{(1^5)}(0)=\varphi_{(1^3)(1^2)}^{(2,1^3)}(0)=\varphi_{(1^3)(1^2)}^{(2^2,1)}(0)=1. $$

\textbf{Step 2.} Note that
$$ \begin{array}{c} F(M(1^3))F(M(2))=u_{(1^3)}(u_{(1^2)}+u_{(2)})=\\
u_{(1^5)}+u_{(2,1^3)}+u_{(2^2,1)}+\varphi_{(1^3)(2)}^{(2,1^3)}(0)u_{(2,1^3)}+\varphi_{(1^3)(2)}^{(3,1^2)}(0)u_{(3,1^2)}.
\end{array} $$
On the other hand
$$F(M(1^3)\ast M((2)))=F(M(3,1^2))= u_{(1^5)}+u_{(2,1^3)}+u_{(2^2,1)}+u_{(3,1^2)}.$$
Therefore 
 $$ \varphi_{(1^3)(2)}^{(3,1^2)}(0)=1 \;\mbox{ and }\; \varphi_{(1^3)(2)}^{(2,1^3)}(0)=0. $$

\textbf{Step 3.} We have
$$ \begin{array}{c} F(M(2,1))F(M(1^2))=(u_{(2,1)}+u_{(1^3)})u_{(1^2)}=\\ 
=\varphi_{(2,1)(1^2)}^{(2,1^3)}(0)u_{(2,1^3)}+\varphi_{(2,1)(1^2)}^{(3,1^2)}(0)u_{(3,1^2)}
+\varphi_{(2,1)(1^2)}^{(2^2,1)}(0)u_{(2^2,1)}+\varphi_{(2,1)(1^2)}^{(3,2)}(0)u_{(3,2)}+\\+u_{(1^5)}+u_{(2,1^3)}+u_{(2^2,1)}.
\end{array} $$
and
$$F(M(2,1)\ast M((1^2)))=F(M(3,2))= u_{(1^5)}+u_{(2,1^3)}+u_{(2^2,1)}+u_{(3,1^2)}+u_{(3,2)}.$$
Therefore 
 $$ \begin{array}{c}
     \varphi_{(2,1)(1^2)}^{(2,1^3)}(0)=\varphi_{(2,1)(1^2)}^{(2^2,1)}(0)=0 \\
     \varphi_{(2,1)(1^2)}^{(3,1^2)}(0)=\varphi_{(2,1)(1^2)}^{(3,2)}(0)=1.
    \end{array}
 $$

\textbf{Step 4.} Finally
$$ \begin{array}{c} F(M(2,1))F(M(2))=(u_{(2,1)}+u_{(1^3)})(u_{(1^2)}+u_{(2)})=\\ 
= \varphi_{(2,1)(2)}^{(2^2,1)}(0)u_{(2^2,1)}+\varphi_{(2,1)(2)}^{(3,2)}(0)u_{(3,2)}
+\varphi_{(2,1)(2)}^{(4,1)}(0)u_{(4,1)}+\varphi_{(2,1)(2)}^{(3,1^2)}(0)u_{(3,1^2)}+
\\+u_{(3,1^2)}+u_{(3,2)}+u_{(3,1^2)}+u_{(1^5)}+u_{(2,1^3)}+u_{(2^2,1)}.
\end{array} $$
and
$$F(M(2,1)\ast M((2)))=F(M(4,1))= u_{(1^5)}+u_{(2,1^3)}+u_{(2^2,1)}+u_{(3,1^2)}+u_{(3,2)}+u_{(4,1)}.$$
Therefore 
 $$ \begin{array}{c}
     \varphi_{(2,1)(2)}^{(2^2,1)}(0)=\varphi_{(2,1)(2)}^{(3,2)}(0)=0 \\
     \varphi_{(2,1)(2)}^{(4,1)}(0)=1 \\
     \varphi_{(2,1)(2)}^{(3,1^2)}(0)=-1.
    \end{array}
 $$
}\end{ex}

\begin{rem}{\rm
 In a~similar way (exploring isomorphism analogous to $F$ given in \cite{wolf})
one may calculate constant terms of Hall polynomials for Dynkin quivers.} 
\end{rem}

\end{document}